\def\draft{0}  

\documentclass[11pt]{article}
\usepackage[letterpaper,hmargin=1in,vmargin=1.25in]{geometry}

\usepackage{amsfonts,url,ifthen,epsfig,eufrak,amsmath,amssymb}

\ifnum\draft=1
\newcommand{\Rnote}[1]{{\bf [Ronen's Note: #1]}}
\newcommand{\Anote}[1]{{\bf [Amir's Note: #1]}}
\else
\newcommand{\Rnote}[1]{}
\newcommand{\Anote}[1]{}
\fi

\newtheorem{theorem}{Theorem}[section]
\newtheorem{definition}[theorem]{Definition}
\newtheorem{lemma}[theorem]{Lemma}

\newtheorem{corollary}[theorem]{Corollary}

\newtheorem{remk}[theorem]{Remark}
\newtheorem{examp}[theorem]{Example}

\def\FullBox{\hbox{\vrule width 8pt height 8pt depth 0pt}}

\def\qed{\ifmmode\qquad\FullBox\else{\unskip\nobreak\hfil
\penalty50\hskip1em\null\nobreak\hfil\FullBox
\parfillskip=0pt\finalhyphendemerits=0\endgraf}\fi}

\def\qedsketch{\ifmmode\Box\else{\unskip\nobreak\hfil
\penalty50\hskip1em\null\nobreak\hfil$\Box$
\parfillskip=0pt\finalhyphendemerits=0\endgraf}\fi}

\newenvironment{proof}{\begin{trivlist} \item {\bf Proof:~~}}
  {\qed\end{trivlist}}

\newcommand{\N}{\mathbb N}

\newcommand{\R}{\mathbb R}

\renewcommand{\d}{\textrm{deg}}

\renewcommand{\Pr}{\mathop{\mathbb{P}}}

\newcommand{\abs}[1]{\left\vert#1\right\vert}

\renewcommand{\abs}[1]{|#1|}

\renewcommand{\d}{\textrm{d}}

\newcommand{\set}[1]{\left\{#1\right\}}
\newcommand{\sett}[2]{ \left\{  #1 \ : \ #2  \right\} }

\newcommand{\E}{\mathop \mathbb{E}}

\title{$t$-Wise Independence with Local Dependencies
\ifnum\draft=1{\\ \small \sc Working Draft, Please Do Not Distribute }\fi }

\author{
Ronen Gradwohl \thanks{Department of Computer Science and Applied
Mathematics, The Weizmann Institute of Science, Rehovot, 76100
Israel. E-mail: \texttt{ronen.gradwohl@weizmann.ac.il}. Research
supported by US-Israel Binational Science Foundation Grant
2002246.} \and
Amir Yehudayoff~\thanks{Department of Computer Science and Applied
Mathematics, The Weizmann Institute of Science, Rehovot, 76100
Israel. E-mail: \texttt{amir.yehudayoff@weizmann.ac.il}.
Research supported by a grant from the Israel Ministry of Science (IMOS) - Eshkol
Fellowship.}
}

\date{}

\begin{document}

\maketitle \thispagestyle{empty}

\begin{abstract}
In this note we prove a large deviation bound on the sum of random
variables with the following dependency structure: there is a
dependency graph $G$ with a bounded chromatic number, in which
each vertex represents a random variable. Variables that are
represented by neighboring vertices may be arbitrarily dependent,
but collections of variables that form an independent set in $G$
are $t$-wise independent.
\end{abstract}


\section{Introduction}
It is often useful to consider a random variable $X= \sum_{i=1}^n X_i$ and bound the
probability that such a sum deviates from its expectation. For independent $X_i$'s,
famous bounds are those of Chernoff \cite{C52} and Hoeffding \cite{H63}.

Sums of variables that are not fully independent but have some sort of a dependency structure
have also been studied -- see for example
Pemmaraju \cite{P01} or the survey of Janson and Ruci\'nski \cite{JR02}. Here we are interested in a
setting studied by Janson \cite{J04}: roughly, in his
formulation there is a dependency graph $G$ of $n$ vertices, in which each variable is represented
by a vertex. Two variables whose corresponding vertices are connected by an edge may be dependent,
whereas independent sets of the graph are independent (see Section~\ref{sec:definitions} for a more
formal description). Janson provides several applications for his bound, such as U-statistics and
the existence of long patterns in random strings (see \cite{J04} for more details).

A different form of dependency structure is motivated by the
computer science literature. This is the setting
in which the $X_i$'s are $t$-wise independent (see for example Bellare and Rompel \cite{BR94}).
This means that every set of $t$
variables is independent, but any $t+1$ may not be.

In this note we consider the situation in which the random variables $X_i$ are dependent
in both fashions: on the one hand, their dependencies are described by a dependency graph $G$.
On the other hand, variables represented by independent sets of the graph are not fully independent,
but only $t$-wise independent. We combine standard techniques used in tail bounds
for sums of $t$-wise independent random variables with the technique of Janson \cite{J04},
and in this manner obtain a tail bound for sums of random variables that are both
$t$-wise independent and have local dependencies.

\subsection{Motivation}

The tail bounds we prove seem to be applicable in numerous
situations, and we will now state two such examples. First, our
inequality was used in a recent game theoretic work of Gradwohl
and Reingold \cite{GR07}. Second,  the inequality can be used in
the hidden pattern problem as a bound on the number of patterns in
a random string that is not fully independent.

\subsubsection{Game theory} In a recent work in game theory, Kalai
\cite{K04} showed that in certain types of large Bayesian games,
the Nash equilibria are not affected by such details as order of
play, possibility of revision, and more. One of the necessary
assumptions in obtaining this result is the independence of
certain random variables related to the players. This was
necessary because of the repeated application of a Chernoff bound.

In a generalization of this work, Gradwohl and Reingold
\cite{GR07} showed how to replace this assumption of independence
by a more general one of limited correlation. One of the tools
used by \cite{GR07} is Corollary~\ref{cor: ber and d}.

\subsubsection{Hidden Pattern Problem} In the hidden pattern
problem, we are given a sequence of $n$ random letters from a
finite alphabet $A$, say $X_1,\ldots,X_n$. Given a word of fixed
length $d$, say $w\in A^d$, one seeks the number of subsequences
$i_1<\ldots<i_d$ such that $X_{i_1}\circ\ldots\circ X_{i_d} = w$.
The case in which the $X_i$'s are independent was studied by
Flajolet et al. as well as Janson \cite{J04}. Bourdon and Vall\'ee
generalize the work by considering strings $X_1,\ldots,X_n$ in
which the $X_i$'s are not fully independent, but rather are
generated by dynamical sources (see \cite{BV02}).
Theorem~\ref{thm: X is concent} can be used in a straightforward
manner to obtain a result similar to that of \cite{J04}, but
applied to variables that are $t$-wise independent.

\section{Definitions}
\label{sec:definitions}

Let $n \in \N$ be an integer. We denote $[n] = \set{1,\ldots,n}$.
For a graph $G$, we denote by $V(G)$ the vertex set of $G$, and by
$E(G)$ the edge set of $G$ (we will consider only simple
undirected graphs). Let $G$ be a graph of size $\abs{V(G)} = n$.
We usually think of $V(G)$ as $[n]$. The following three
definitions are standard graph definitions.

\begin{definition}[independent set]
$S \subseteq V(G)$ is an \emph{independent set} of vertices in $G$
if no two vertices in $S$ share an edge (according to $G$).
\end{definition}

\begin{definition}[coloring]
For $k \in \N$, a \emph{$k$-coloring} of $G$ is a map from
$V(G)$ to $[k]$ such that each two adjacent vertices are mapped to
different integers.
\end{definition}

\begin{definition}[chromatic number]
The \emph{chromatic number} of $G$, denoted by $\chi(G)$, is the smallest integer $k$ such that there
exists a $k$-coloring of $G$.
\end{definition}

Note that if the degree of $G$ is at most $d \in \N$, then
$\chi(G) \leq d+1$ (since the greedy algorithm for coloring works).

Next, we give three definitions concerning distributions: the
first one is the standard definition of $t$-wise independence, the
second definition is of a dependency graph, and the third
definition, which combines the first two definitions, is of the
family of distributions for which our tail bounds apply.

\begin{definition}[$t$-wise independence]
For $m,t \in \N$, the random variables $Y_1,\ldots,Y_m$ are
\emph{$t$-wise independent}, if for every $T \subseteq [m]$ of
size $t$ the set of variables $\sett{Y_i}{i \in T}$ is
independent.
\end{definition}

\begin{definition}[agree]
Let $n \in \N$, and let $G$ be a graph of size $n$. We say that
the random variables $X_1,\ldots,X_n$ \emph{agree} with the graph
$G$, if for every independent set of vertices $S \subseteq V(G)$,
the set of variables $\sett{X_i}{i \in S}$ is independent ($G$ is
sometimes called a dependency graph).
\end{definition}

\begin{definition}[$t$-agree]
For a dependency graph $G$ as above, we say that the random
variables $X_1,\ldots,X_n$ \emph{$t$-agree} with $G$, if for every
independent set of vertices $S \subseteq V(G)$, the set of
variables $\sett{X_i}{i \in S}$ is $t$-wise independent.
\end{definition}

\section{Results}
We are now ready to state our main result, which is a large
deviation bound on the sum of random variables that $t$-agree with
a graph $G$ of chromatic number $\chi(G)$.
\begin{theorem}
\label{thm: X is concent} Let $n,t \in \N$ be such that $t >0$ is
even. Let $G$ be a graph of size $n$, and let $X_1,\ldots,X_n$ be
random variables that take values in $[0,1]$ and $t$-agree with
$G$. Let $X = \sum_{i \in [n]} X_i$ and let $\mu = \E[X]$. Then,
for every positive real $a > 0$,
$$\Pr\left[\abs{X-\mu} \geq a\right]  < 2 \sqrt{\pi t} \cdot \left( \frac{
  \sqrt{n t\cdot\chi(G)}}{a}
 \right)^t.$$
\end{theorem}

\noindent When the random variables are Bernoulli and the graph is of
bounded degree we have the following corollary:

\begin{corollary}
\label{cor: ber and d} Let $n,d,t \in \N$ be such that $t >0$ is
even. Let $G$ be a graph of size $n$ and degree at most $d$. Let
$p \in (0,1)$, and let $X_1,\ldots,X_n$ be $Be(p)$ random
variables that $t$-agree with $G$. Let $X = \sum_{i \in [n]} X_i$.
Then for every positive real $a > 0$,
$$\Pr\left[X \geq (1+a) pn\right]  < 2 \sqrt{\pi t} \cdot \left( \frac{
 \sqrt{(d+1) \cdot t}}{ap \sqrt{n}}
 \right)^t.$$
\end{corollary}
The same bound holds for $\Pr\left[X \leq (1-a) pn\right]$.

\section{Proof of Main Result}
In our proof we will need to bound the $t$-moment of the sum of $t$-wise independent
random variables. The following bound is well known -- see Bellare and Rompel \cite{BR94}
for a proof.

\begin{lemma}
\label{lem: exp x to k is small} Let $m,t \in \N$ be such that $t
> 0$ is even, and let $Y_1,\ldots,Y_m$ be $t$-wise independent
random variables taking values in $[0,1]$. Let $Y = \sum_{i \in
[m]} Y_i$ and let $\mu = \E\left[Y\right]$. Then
$$\E\left[(Y-\mu)^t\right] < 2 \cdot e^{\frac{1}{6t}} \cdot \sqrt{\pi t}
\cdot \left( \frac{mt}{e} \right)^{t/2}.$$
\end{lemma}

\noindent We now prove Theorem~\ref{thm: X is concent}.

\begin{proof}
Let $G$ be a graph of size $n$, and let $X_1,\ldots,X_n$ be random
variables that $t$-agree with $G$. Let $X = \sum_{i \in [n]} X_i$
and let $\mu = \E\left[X\right]$. Let $f$ be a $k$-coloring of $G$ such that
$k = \chi(G)$. For every $j \in [k]$, denote
$$V_j = f^{-1}(j), $$ which is an independent set of vertices. So
for all $j \in [k]$, the set of variables $\sett{X_i}{i \in V_j}$
is $t$-wise independent. By Lemma~\ref{lem: exp x to k is
small}, for every $j \in [k]$,
$$\E\left[(Y_j-\mu_j)^t\right] < 2 \cdot e^{\frac{1}{6t}} \cdot \sqrt{\pi t}
\cdot \left( \frac{\abs{V_j} t}{e} \right)^{t/2},$$ where $Y_j =
\sum_{i \in V_j} X_i$ and $\mu_j = \E\left[Y_j\right]$.

We now bound the $t$-moment of $X$. Let $p_1,\ldots,p_k$ be $k$
non-negative real numbers such that $\sum_{j \in [k]} p_j = 1$ (to
be determined later). By Jensen's inequality and linearity
of expectation,
\begin{eqnarray*} \E\left[(X - \mu)^t\right] & = & \E\left[\left(\sum_{j \in [k]} p_j
\frac{Y_j - \mu_j}{p_j}\right)^t\right]
\\ & \leq & \sum_{j \in [k]} p_j \frac{\E\left[(Y_j-\mu_j)^t\right]}{p_j^t} \\ &
< & \sum_{j \in [k]} p_j \frac{2 \cdot e^{\frac{1}{6t}} \cdot
\sqrt{\pi t} \cdot \left( \frac{\abs{V_j} t}{e}
\right)^{t/2}}{p_j^t}.
\end{eqnarray*}
For $j \in [k]$, set $q_j \in \R$ to be such that $$q_j^t = 2
\cdot e^{\frac{1}{6t}} \cdot \sqrt{\pi t} \cdot \left(
\frac{\abs{V_j} t}{e} \right)^{t/2},$$ and set $$p_j =
\frac{q_j}{\sum_{\ell \in [k]} q_\ell}.$$

Substituting these values of $p_j$ yields
$$\E\left[(X - \mu)^t\right] < \left( \sum_{j \in [k]} q_j \right)^t =
\left(\sum_{j \in [k]} (2 \cdot e^{\frac{1}{6t}} \cdot \sqrt{\pi
t})^{1/t} \cdot \left( \frac{\abs{V_j} t}{e} \right)^{1/2}
\right)^t,$$ which implies that
$$\E\left[(X - \mu)^t\right] < 2 \cdot \sqrt{\pi
t} \cdot \left(\sqrt{k n t} \right)^t$$ by Cauchy-Schwartz. Since $t$ is even we now use
Markov's Inequality, implying that for every real number $a
> 0$,
$$\Pr\left[\abs{X-\mu} \geq a\right] = \Pr\left[(X-\mu)^t \geq a^t\right] \leq
\frac{\E\left[(X-\mu)^t\right]}{a^t} < 2 \sqrt{\pi t} \cdot \left( \frac{\sqrt{k n t}}{a}
 \right)^t.$$ Substituting $k = \chi(G)$, we get that
$$\Pr\left[\abs{X-\mu} \geq a\right] < 2 \sqrt{\pi t} \cdot \left( \frac{
 \sqrt{ n t\cdot\chi(G)}}{a}
 \right)^t.$$
\end{proof}


\begin{thebibliography}{99}

\bibitem{BR94} M. Bellare and J. Rompel.
Randomness-efficient oblivious sampling.
Proceedings 35th Annual Symposium on the Foundations of Computer Science, IEEE, 1994.

\bibitem{BV02} J. Bourdon and B. Vall\'ee. Generalized pattern matching statistics.
In Mathematics and Computer Science (Colloquium Proceedings, Versailles, 2002).
B. Chauvin et al. editors, Birkh\"{a}user Verlag, 2002: 229-245.

\bibitem{C52} H. Chernoff.
A measure of asymptotic efficiency for tesets of a hypothersis based on the sum
of observations.
Annals of Mathematical Statistics 23 (1952): 493-507.

\bibitem{FSV06} P. Flajolet, W. Szpankowski, and B. Vall\'ee. Hidden word statistics.
In Journal of the ACM, Volume 53:1, January 2006, pages 147--183.

\bibitem{GR07} R. Gradwohl and O. Reingold.
Partial exposure and correlated types in large games. Submitted

\bibitem{H63} W. Hoeffding.
Probability inequalities for sums of bounded random variables.
Journal of the American Statistical Association 58 (1963): 13-30.

\bibitem{J04}
S. Janson. Large deviations for sums of partly dependent random variables.
Random Structures and Algorithms 24 (2004): 234--248.

\bibitem{JR02} S. Janson and A. Ruci\'nski. The infamous upper tail.
Random Structures and Algorithms 20 (2002): 317-342.

\bibitem{K04} E. Kalai. Large robust games.
Econometrica, Vol. 72, No. 6, November 2004.
Pages 1631-1665.

\bibitem{P01} S.V. Pemmaraju. Equitable colong extends Chernoff-Hoeffding bounds.
Approximation, Randomization and Combinatorial Optimization: Algorithms and Techniques
(APPROX 2001 and RANDOM 2001), eds. M.X. Goemans et al.
Lecture Notes in Computer Science 2129, Springer, 2001: 285-296.

\end{thebibliography}
\end {document}